\documentclass[11pt, twoside]{article}
\usepackage{latexsym}
\usepackage{amsmath}
\usepackage{amssymb}
\usepackage[all]{xy}
\usepackage{amsfonts}
\usepackage{verbatim}
\usepackage{amsthm}
\usepackage{mathrsfs}
\usepackage{epsfig}
\usepackage{xy}
\usepackage{array}
\usepackage{stmaryrd}
\usepackage{graphicx,color}
\usepackage{xcolor}
\usepackage[colorlinks=true,linkcolor=blue,citecolor=blue]{hyperref}
\usepackage{tikz}
\usetikzlibrary{arrows,calc}
\usepackage{etex}
\usepackage{mathdots}
\usepackage{float}
\usepackage{graphics}
\usepackage{pdflscape}
\usepackage{bm}
\usepackage{anysize,hyperref}
\input xypic
\xyoption{all}

\usepackage[perpage,symbol]{footmisc}
\usepackage{setspace}
\SelectTips{cm}{10}

\begin{document}
\title{\Large{\bf On Nichols bicharacter algebras$^\bigstar$\footnotetext{
\hspace{-1em}$^\bigstar$Weicai Wu was supported by the Hunan Provincial Natural Science Foundation of China (Grant No. 2020JJ5210).}}}
\medskip
\author{Weicai Wu}

\date{}

\maketitle
\def\blue{\color{blue}}
\def\red{\color{red}}

\newtheorem{Proposition}{Proposition}[section]
\newtheorem{Theorem}[Proposition]{Theorem}
\newtheorem{Definition}[Proposition]{Definition}
\newtheorem{Corollary}[Proposition]{Corollary}
\newtheorem{Lemma}[Proposition]{Lemma}
\newtheorem{Example}[Proposition]{Example}
\newtheorem{Remark}[Proposition]{Remark}
\newtheorem{Conjecture}[Proposition]{Conjecture}
\newtheorem{Question}[Proposition]{Question}

\baselineskip=17pt
\parindent=0.5cm
\vspace{-6mm}

\begin{abstract}
\baselineskip=16pt
In this paper we define two Lie operations, and with that we define the bicharacter algebras, Nichols bicharacter algebras, etc. We obtain explicit bases for $\mathfrak L(V)${\tiny $_{R}$} and $\mathfrak L(V)${\tiny $_{L}$} over a connected braided vector $V$ of  diagonal type with $\dim V=2$ and $p_{1,1}=p_{2,2}= -1$. We give the sufficient and necessary conditions for  $\mathfrak L(V)${\tiny $_{R}$}$= \mathfrak L(V)$, $\mathfrak L(V)${\tiny $_{L}$}$= \mathfrak L(V)$, $\mathfrak B(V) = F\oplus \mathfrak L(V)${\tiny $_{R}$} and $\mathfrak B(V) = F\oplus \mathfrak L(V)${\tiny $_{L}$}, respectively.
We show that if $\mathfrak B(V)$ is a connected Nichols algebra of diagonal type with $\dim V>1$, then $\mathfrak B(V)$ is finite-dimensional if and only if $\mathfrak L(V)${\tiny $_{L}$} is finite-dimensional if and only if $\mathfrak L(V)${\tiny $_{R}$} is finite-dimensional.\\[0.2cm]
\end{abstract}

\pagestyle{myheadings}
\markboth{\rightline {\scriptsize   Weicai Wu}}
         {\leftline{\scriptsize On Nichols bicharacter algebras}}

\section{Introduction}\label{s0}
Nichols algebras play a central role in the theory of pointed Hopf algebras, and
the classification of (braided) Lie algebras has been an international hotspot and cutting-edge problem in the field of algebras. A great deal of attention has been paid to the question of finite-dimensionality of Nichols algebras (see \cite {AHS08,AS10,He05,He06a,He06b,WZZ15b}).
The interest in this problem arose from the work of Andruskiewitsch and Schneider \cite{AS98} on classification of finite dimensional (Gelfand-Kirillov) pointed Hopf algebras which are generalizations of quantized  enveloping algebras of semi-simple Lie algebras. On the other hand, Lie algebra arising from a Nichols algebra was studied in \cite {AAB17, AAB19}.
The theory of Lie superalgebras has been developed systematically, which includes the
representation theory and classifications of simple Lie superalgebras and their varieties
\cite{Ka77}. A sophisticated
multilinear version of the Lie bracket was considered in \cite {{Kh99, Pa98, LR95}}. Various generalized
Lie algebras have already appeared under different names, e.g. Lie color algebras, $\epsilon $ Lie
algebras \cite{Sc79}, quantum and braided Lie algebras \cite {Ma94}, braided $m$-Lie algebras  \cite{ZZ04} and
generalized Lie algebras \cite{BFM96}.

In this work we define two new Lie operations, this is one of the main results in this paper, which enables us to define the bicharacter algebras(all linear Lie algebras and Kac-Moody algebras are bicharacter algebras), Nichols bicharacter algebras, etc. Braided Lie algebra (with braided Lie operation $[u, v]=vu-p_{v,u}uv$) seems to be `located' between two Lie algebras(the usual ones with Lie operation $[u, v]^{-}=uv -vu$ and bicharacter algebras with Lie operation $[u, v]${\tiny $_{L}$} $= p_{v,u}uv -p_{u,v}vu$ or $[u, v]${\tiny $_{R}$} $= p_{u,v}uv-p_{v,u}vu$), which is of great significance for revealing the relationship between braided Lie algebra and Lie algebra. In the remaining part of this paper will focus on Nichols bicharacter algebras, which is a cross research targe in Lie algebras and Nichols algebras.
It was proven that a Nichols algebra is finite-dimensional if and only if the corresponding Nichols bicharacter algebras is finite-dimensional.
This provides a new method for determining when a Nichols algebra is finite dimensional.

This paper is organized as follows. In Section \ref{s1} we recall some results on Lie algebras and Nichols algebras, and define bicharacter algebras and Nichols bicharacter algebras. Section \ref{s3} presents explicit bases for $\mathfrak L(V)${\tiny $_{L}$} and
$\mathfrak L(V)${\tiny $_{R}$},  where  $V$ is a connected braided vector space of  diagonal type with $\dim V=2$ and $p_{1,1}=p_{2,2}= -1$.
In Section \ref{s4} we present the sufficient and necessary conditions for  $\mathfrak L(V)${\tiny $_{R}$}$= \mathfrak L(V)$, $\mathfrak L(V)${\tiny $_{L}$}$= \mathfrak L(V)$, $\mathfrak B(V) = F\oplus \mathfrak L(V)${\tiny $_{R}$} and $\mathfrak B(V) = F\oplus \mathfrak L(V)${\tiny $_{L}$}, respectively.
In Section \ref{s5} we prove that $\mathfrak B(V)$ is finite-dimensional if and only if $\mathfrak L(V)${\tiny $_{L}$} is finite-dimensional
if and only if $\mathfrak L(V)${\tiny $_{R}$} is finite-dimensional, when $\mathfrak B(V)$ is a connected Nichols algebra of diagonal type with $\dim V>1$.\\[1mm]
{\bf Conventions.}  $\mathbb Z =: \{x \mid x \hbox { is an  integer}\}$.
$\mathbb N_0 =: \{x \mid x \in \mathbb Z,  x\ge 0\}.$
$\mathbb N =: \{x \mid x \in \mathbb Z,  x>0\}.$  $F$ denotes the base field of characteristic zero.
$F^{*}=F\backslash\{0\}$.

\section{Preliminaries}\label {s1}
In this section we recall some results on Lie algebras and Nichols algebras, and define bicharacter algebras and Nichols bicharacter algebras.

Let $E=\{e_{1}, e_{2}, \cdots, e_{n}\}, e_i =: ( \stackrel{i} { \overbrace{0,  \cdots, 0,  1}},  \cdots, 0)\in \mathbb Z^{n}$ ,  $1\le i \le n.$
Let $A := \{ x_{1}, x_{2}, \cdots, x_{n} \}$  be an alphabet and a basis of $V = {\rm span } (A)$, likewise, associated with each letter $x_{i}$ are an element $e_{i}=\deg(x_{i})$ of $\mathbb Z^{n}$ and a character $\chi^{x_{i}} : \mathbb Z^{n} \longrightarrow F^{*}$. Let an empty word by $1$, $W:= \{u \mid u \hbox { is a word} \}$ and a basis of $T(V) = {\rm span } (W)$. Let $V^{*}$ denote the vector space dual to $V$ and $\{ y_{1}, y_{2}, \cdots, y_{n} \}$ denote the dual basis of $V^{*}$, where $\deg(y_{i})=-e_{i},\chi^{y_{i}} : \mathbb Z^{n} \longrightarrow F^{*}$. Let $T(\overline{V})$ is free generated by $1\cup \overline{V}:=\{1, x_{i}, y_{i}\mid 1\leq i\leq n\}$ and $\overline{W}$ is a basis of $T(\overline{V})$ with that each element in $\overline{W}$ is gotten from one in $1\cup \overline{V}$ by multiplication elements of $1\cup \overline{V}$. $\overline{W}$ is a monoid with $1$, for all $u\in \overline{W}$, denote by $\deg(u)$ an element of the group $\mathbb Z^{n}$ which results from $u$ by replacing each occurrence of the letter $x_{i}$ (or $y_{i}$) with $e_{i}$ (or $-e_{i}$). By $\chi^{u}$ we denote a character which results from $u$ by replacing all $x_i$ (or $y_{i}$) with $\chi^{x_{i}}$ (or $\chi^{y_{i}}$) . For a pair of words $u,v\in \overline{W}$, put $p_{u,v}=\chi^{u}(\deg{v})$, where $p_{x_{i},x_{j}}=\chi^{x_{i}}(\deg{x_{j}}),p_{x_{i},x_{j}}^{-1}=\chi^{x_{i}}
(\deg{y_{j}}),
p_{x_{i},x_{j}}^{-1}=\chi^{y_{i}}(\deg{x_{j}}),p_{x_{i},x_{j}}
=\chi^{y_{i}}(\deg{y_{j}})$. Obviously, $p_{uu_1, v} = p_{u,v}p_{u_{1}, v},p_{u, vv_{1}} = p_{u,v}p_{u,v_{1}}$, that is the operator $p_{\cdot,\cdot}$ is a bicharacter defined on $\overline{W}$. Moreover, we call $p_{\cdot,\cdot}$ is a symmetrical bicharacter if $p_{u,v}=p_{v,u}$ for all $u,v\in \overline{W}$.
Define
$u$ $\circ${\tiny $ _{L}$}$v:=p_{v,u}uv$, $u$ $\circ${\tiny $ _{R}$}$v:=p_{u,v}uv$. It is clear ($T(\overline{V})$, $\circ${\tiny $_{L}$}) and ($T(\overline{V})$, $\circ${\tiny $_{R}$}) are associative algebras. Assume that $A(V)$ is in general a quotient algebra of some subalgebras of $T(\overline{V})$, then ($A(V)$, $\circ${\tiny $_{L}$}) and ($A(V)$, $\circ${\tiny $_{R}$}) are associative algebras. Especially, ($T(V)$, $\circ${\tiny $_{L}$}) and ($T(V)$, $\circ${\tiny $_{R}$}) are associative algebras. Let $I$ be the (two sided) ideal in $T(V)$ generated by all $x_{i}x_{j}-x_{j}x_{i},1\leq i,j\leq n$ and call $S(V)=T(V)/I$ the symmetric algebra on $V$.

\begin{Lemma} \label {0} If  $(A,\circ)$ is an associative algebra together with the bilinear operation $[\cdot , \cdot]$ defined by $[u, v]= u\circ v -v\circ u$ for any homogeneous elements $u, v \in A$, then $(A, [\cdot,\cdot])$ is a Lie algebra.
\end{Lemma}

If we define the bilinear operations $[u, v]${\tiny $_{L}$}$=u$ $\circ${\tiny $ _{L}$}$v-v$ $\circ${\tiny $ _{L}$}$u$ and $[u, v]${\tiny $_{R}$} $=u$ $\circ${\tiny $ _{R}$}$v-v$ $\circ${\tiny $ _{R}$}$u$ for all $u,v\in W$, then ($A(V)$, $[\cdot , \cdot]${\tiny $_{L}$}) and ($A(V)$, $[\cdot , \cdot]${\tiny $_{R}$}) are Lie algebras, which are called L-bicharacter algebra and R-bicharacter algebra, written as $A(V)${\tiny $_{L}$} and $A(V)${\tiny $_{R}$}, respectively. Both L-bicharacter algebra and R-bicharacter algebra are called bicharacter algebras.

\begin{Remark} \label{1.1} {\rm (i)} ($S(V)$, $\circ${\tiny $_{L}$})(resp. ($S(V)$, $\circ${\tiny $_{R}$})) is a commutative associative algebra if and only if $p_{\cdot,\cdot}$ is a symmetrical bicharacter.

{\rm (ii)} Let  $S(V)${\tiny $_{L}$} and $S(V)${\tiny $_{R}$} denote the  Lie algebras generated by $V$ in $S(V)$ under Lie operations $[u, v]${\tiny $_{L}$} $= p_{v,u}uv -p_{u,v}vu$, $[u, v]${\tiny $_{R}$} $= p_{u,v}uv-p_{v,u}vu$, respectively, for any homogeneous elements $u, v\in S(V)$. Then $[u, v]${\tiny $_{L}$}$= (p_{v,u}-p_{u,v})uv=-[u, v]${\tiny $_{R}$} and  $S(V)${\tiny $_{L}$} $=S(V)${\tiny $_{R}$}.

{\rm (iii)} Linear Lie algebras and Kac-Moody algebras are bicharacter algebras, where $p_{v,u}= p_{u,v}=1$.
\end{Remark}

Let $h_{i}=:x_{i}y_{i}-y_{i}x_{i}$, $<\alpha_{j},\alpha_{i}>$ are Cartan integers, $1\leq i,j\leq n$. $(ad \ u)${\tiny $_{L}$}$(v)=:[u, v]${\tiny $_{L}$}, $(ad\ u)${\tiny $_{R}$}$(v)=:[u, v]${\tiny $_{R}$}, for any $u, v\in T(\overline{V})$.

\begin {Proposition} \label{1.1''} {\rm (See \cite [Prop. 18.1, Theor. 18.3]{Hu72})} {\rm (1)} With the above notation, $(A_{V})${\tiny $_{L}$} is generated by
$\{x_{i},y_{i}|1\leq i\leq n\}$ and $p_{\cdot,\cdot}$ is a symmetrical bicharacter, then $(A_{V})${\tiny $_{L}$} is a finite dimensional semisimple Lie algebra if and only if these generators satisfy the following relations:

{\rm (i)} $[h_{i},h_{j}]${\tiny $_{L}$}$=0$, $1\leq i,j\leq n$.

{\rm (ii)} $[x_{i},y_{i}]${\tiny $_{L}$}$=p_{x_{i},x_{i}}^{-1}h_{i}$, $[x_{i},y_{j}]${\tiny $_{L}$}$=0$ if $i\neq j$.

{\rm (iii)} $[h_{i},x_{j}]${\tiny $_{L}$}$=<\alpha_{j},\alpha_{i}>x_{j}$,
$[h_{i},y_{j}]${\tiny $_{L}$}$=-<\alpha_{j},\alpha_{i}>y_{j}$, $1\leq i,j\leq n$.

{\rm (iv)} $(ad \ x_{i})${\tiny $_{L}$}$^{-<\alpha_{j},\alpha_{i}>+1}(x_{j})=0$ if $i\neq j$.

{\rm (v)} $(ad \ y_{i})${\tiny $_{L}$}$^{-<\alpha_{j},\alpha_{i}>+1}(y_{j})=0$ if $i\neq j$.

{\rm (2)} With the above notation, $(A_{V})${\tiny $_{R}$} is generated by
$\{x_{i},y_{i}|1\leq i\leq n\}$ and $p_{\cdot,\cdot}$ is a symmetrical bicharacter, then $(A_{V})${\tiny $_{R}$} is a finite dimensional semisimple Lie algebra if and only if these generators satisfy the following relations:

{\rm (i)} $[h_{i},h_{j}]${\tiny $_{R}$}$=0$, $1\leq i,j\leq n$.

{\rm (ii)} $[x_{i},y_{i}]${\tiny $_{R}$}$=p_{x_{i},x_{i}}^{-1}h_{i}$, $[x_{i},y_{j}]${\tiny $_{R}$}$=0$ if $i\neq j$.

{\rm (iii)} $[h_{i},x_{j}]${\tiny $_{R}$}$=<\alpha_{j},\alpha_{i}>x_{j}$,
$[h_{i},y_{j}]${\tiny $_{R}$}$=-<\alpha_{j},\alpha_{i}>y_{j}$, $1\leq i,j\leq n$.

{\rm (iv)} $(ad \ x_{i})${\tiny $_{R}$}$^{-<\alpha_{j},\alpha_{i}>+1}(x_{j})=0$ if $i\neq j$.

{\rm (v)} $(ad \ y_{i})${\tiny $_{R}$}$^{-<\alpha_{j},\alpha_{i}>+1}(y_{j})=0$ if $i\neq j$.
\end {Proposition}

\begin {Example} \label{1.1''''} When $\dim V=1$, $(A_{V})${\tiny $_{R}$} is a Lie algebra with generators $x_{1},y_{1}$ and relations $[x_{1},y_{1}]${\tiny $_{R}$}$=p_{x_{1},x_{1}}^{-1}h_{1}$, $[h_{1},x_{1}]${\tiny $_{R}$}$=2x_{1}$,
$[h_{1},y_{1}]${\tiny $_{R}$}$=-2y_{1}$. Moreover, $(A_{V})${\tiny $_{R}$}$=(A_{V})${\tiny $_{L}$}.

\end {Example}

From now on let $V=\{x_{1},   \cdots,   x_{n}\}$ is  a basis of  vector space  $V$ and
$C(x_{i}\otimes  x_{j}) =q_{ij} x_{j}\otimes x_{i}$ with $q_{ij} =p_{x_{i},x_{j}}$,
then $V$  is called a braided vector space of diagonal type, $\{x_1,   \cdots,   x_n\}$
is called canonical basis and $(q_{ij})_{n\times n}$  is called braided matrix. In this case, $x_{i}$ $\circ${\tiny $ _{L}$}$x_{j}=({g_{j}}_{\cdot}x_{i})x_{j}$, $x_{i}$ $\circ${\tiny $ _{R}$}$x_{j}=x_{i}({g_{i}}_{\cdot}x_{j})$.
Throughout this paper all of braided vector spaces are connected and  of diagonal type without special announcement.

Let $S_{m}\in End_{k}(V^{\otimes m})$
and $S_{1, j}\in End_{k}(V^{\otimes j+1})$ denote the maps

$S_{m}=\prod \limits_{j=1}^{m-1}(id^{\bigotimes m-j-1}\bigotimes S_{1, j})$,

$S_{1, j}=id+C_{12}^{-1}+C_{12}^{-1}C_{23}^{-1}+\cdots+C_{12}^{-1}C_{23}^{-1}\cdots C_{j, j+1}^{-1}$
(in leg notation) for $m\geq 2$ and $j\in \mathbb N$. Then the subspace
$S=\bigoplus \limits _{m=2}^{\infty} ker S_{m}$ of the tensor algebra $T(V)$
is a two-sided ideal,  and algebra $\mathfrak B(V)=T(V)/S$ is termed the Nichols algebra associated to $(V, C)$ (See \cite [Def. 1.2.2]{He05}).

Let $R_{m} := \{ \alpha \mid \alpha$ is a primitive $m$-th root of $1 \}$. Set
$R_{3} := \{ 1,\omega,\omega^{2} \}$. Let $i$ denote $x_{i}$ in short, sometimes.

\begin{Lemma} \label {1.4} Assume that $i\neq j$.

{\rm (i)}
$[x_{i}, x_{j}]${\tiny $_{L}$} $= 0$ if and only if $p_{i,j}=p_{j,i}=1$.

{\rm (ii)}
$[x_{i}, x_{j}]${\tiny $_{R}$} $= 0$ if and only if $p_{i,j}=p_{j,i}=1$ or $p_{i,j}=\omega,p_{j,i}=\omega^{2}$ or $p_{i,j}=\omega^{2},p_{j,i}=\omega$.
\end{Lemma}
\noindent {\it Proof.}
{\rm (i)} $<y_{i},[x_{i}, x_{j}]${\tiny $_{L}$}$>=(p_{j,i}-1)x_{j}$,
$<y_{j},[x_{i}, x_{j}]${\tiny $_{L}$}$>=(1-p_{i,j})x_{i}$.

{\rm (ii)} $<y_{i},[x_{i}, x_{j}]${\tiny $_{R}$}$>=p_{i,j}^{-1}(p_{i,j}^{2}-p_{j,i})x_{j}$,
$<y_{j},[x_{i}, x_{j}]${\tiny $_{R}$}$>=p_{j,i}^{-1}(p_{i,j}-p_{j,i}^{2})x_{i}$.
\qed

Let  $\mathfrak L^{-}(V)$, $\mathfrak L(V)${\tiny $_{L}$} and $\mathfrak L(V)${\tiny $_{R}$} denote the  Lie algebras generated by $V$ in $\mathfrak B(V)$ under Lie operations $[u, v]^{-}=uv-vu$, $[u, v]${\tiny $_{L}$} $= p_{v,u}uv -p_{u,v}vu$, $[u, v]${\tiny $_{R}$} $= p_{u,v}uv-p_{v,u}vu$, respectively, for any homogeneous elements $u, v\in \mathfrak B(V)$.
$(\mathfrak L^ {-} (V),    [\cdot , \cdot ]^ {-} )$,  ($\mathfrak L(V)${\tiny $_{L}$}, $[\cdot , \cdot]${\tiny $_{L}$})
and
($\mathfrak L(V)${\tiny $_{R}$}, $[\cdot , \cdot]${\tiny $_{R}$}) are called Nichols Lie algebra,  Nichols L-bicharacter algebra
and Nichols R-bicharacter algebra of $V$, respectively.
Let  $\mathfrak L(V)$ denote the  braided Lie algebras generated by $V$ in $\mathfrak B(V)$ under Lie operations $[u, v]=vu-p_{u,v}vu$,  for any homogeneous elements $u, v\in \mathfrak B(V)$. $(\mathfrak L(V), [\cdot , \cdot])$ is called Nichols braided Lie algebra of $V$.

\begin{Remark} \label{1.3} Assume that $p_{\cdot,\cdot}$ is a symmetrical bicharacter. Then $\mathfrak L^ - (V)=\mathfrak L(V)${\tiny $_{L}$} $=\mathfrak L(V)${\tiny $_{R}$} since $[u,v]${\tiny $_{R}$}
$=[u,v]${\tiny $_{L}$}$=p_{u,v}[u,v]^{-}$ for any homogeneous elements $u,v \in \mathfrak B(V)$.
\end{Remark}

Let $l_{u}^0[v]${\tiny $_{L}$}$:= v$,  $l_{u}^{i}[v]${\tiny $_{L}$}$:=  [u, l_{u}^{i-1} [ v]${\tiny $_{L}$}$]${\tiny $_{L}$}, $i\geq 1$. Similarly define $r_{u}^i[v]${\tiny $_{L}$}.
Let $l_{u}^0[v]${\tiny $_{R}$}$:=  v$,  $l_{u}^{i}[v]${\tiny $_{R}$}$:=  [u, l_{u}^{i-1} [ v]${\tiny $_{R}$}$]${\tiny $_{R}$}, $i\geq 1$. Similarly define $r_{u}^{i}[v]${\tiny $_{R}$}. In fact,
$l_{u}^{i}[v]${\tiny $_{R}$}$:= [u, [u, \cdots, [u, v]${\tiny $_{R}$}$\cdots ]${\tiny $_{R}$}$]${\tiny $_{R}$}.

\section{Nichols bicharacter algebra with $\dim V=2$ and $p_{1,1}=p_{2,2}=-1$}\label {s3}
In this section we obtain explicit bases for $\mathfrak L(V)${\tiny $_{L}$} and
$\mathfrak L(V)${\tiny $_{R}$}  where  $V$ is a connected braided vector space of  diagonal type with $\dim V=2$ and $p_{1,1}=p_{2,2}= -1$.

\begin {Lemma}\label {4.1} {\rm (See \cite[Lemma 4.1]{WZZ})}  Assume that $\mathfrak B(V)$ is a connected Nichols algebra of diagonal type with $\dim V=2$ and $p_{1,1}=p_{2,2}= -1$.

{\rm (i)} If $u$ is a non-zero monomial, then there exists $\alpha \in F^{*}$ such that $u =  \alpha (x_1x_2)^k x_1, $ or $\alpha x _2(x_1x_2)^k ,$
 or $\alpha (x_1x_2)^{k+1},$ or $\alpha (x_2x_1)^{k+1},$ $k\ge 0.$

{\rm (ii)} Let  $P :=  \begin{array}{ll}
\{ x_2(x_{1}    x_2)^{k},(x_{1}  x_2)^{k}x_1, (x_{1}    x_2)^{k},  x_{2}( x_1   x_2)^{k}x_1, 0 \le k <  {\rm ord}
(p_{1,2} p_{2,1})\}.
\end{array} $ Then $P$ is a basis of $\mathfrak B(V)$.
\end {Lemma}

\subsection{Nichols R-bicharacter algebra}

\begin {Lemma}\label {4.2}   Assume that $\mathfrak L(V)${\tiny $_{R}$} is a connected Nichols R-bicharacter algebra of diagonal type with $\dim V=2$ and $p_{1,1}=p_{2,2}= -1$.

{\rm (i)} $(r_{x_{1}}r_{x_{2}})^{i}[x_{1}]${\tiny $_{R}$}$ = (-2)^{i} (p_{1,2}p_{2,1})^{\frac{i(i+1)}{2}}(x_{1}x_{2})^{i}x_{1}$ and

$(r_{x_{2}}r_{x_{1}})^{i}[x_{2}]${\tiny $_{R}$}$ = (-2)^{i} (p_{1,2}p_{2,1})^{\frac{i(i+1)}{2}}x_{2}(x_{1}x_{2})^{i}$ for $i \ge 0.$

{\rm (ii)} For any method $\sigma$ adding bracket  $[\cdot,\cdot]${\tiny $_{R}$} and  $k >0$, there exist $\alpha_k, \beta _k \in F$
such that
$\sigma (x_{i_1}, x_{i_2}, \cdots , x_{i_{2k+1}}) = \beta_k  (x_1x_2)^ {k} x_1$ or $ \beta_k  x_2(x_1x_2)^ {k}$
and $\sigma (x_{i_1}, x_{i_2}, \cdots , x_{i_{2k}}) = \alpha_k  ( p_{1,2}^{k}(x_1x_2)^ k - p_{21}^{k}(x_2x_1)^ k)$.

{\rm (iii)} $
[[ [[x_{1}, x_{2}]${\tiny $_{R}$}$, x_1 ]${\tiny $_{R}$}$, x_2]${\tiny $_{R}$}$\cdots , x_{1}]${\tiny $_{R}$}$]${\tiny $_{R}$}$ = 2^{k-1}
 (p_{1,2}p_{2,1})^{\frac{k(k-1)}{2}}(p_{1,2}^{k} (x_{1}x_{2})^{k}- p_{2,1}^{k}(x_{2}x_{1})^{k}   )$,
where there exist $2k-1$ brackets in the left hand side.

{\rm (iv)} $$<y_1, p_{1,2}^{i}(x_{1}x_2)^{i}-p_{2,1}^{i}(x_{2}x_1)^{i}>=p_{1,2}^{i}
(1+(-p_{1,2}^{-2}p_{2,1})^{i})x_2(x_{1}    x_2)^{i-1},$$
$$<y_2, p_{1,2}^{i}(x_{1}x_2)^{i}-p_{2,1}^{i}(x_{2}x_1)^{i}>=-p_{2,1}^{i}
(1+(-p_{1,2}p_{2,1}^{-2})^{i})(x_{1}    x_2)^{i-1}x_1.$$
\end {Lemma}
\noindent {\it Proof.}
{\rm (i)} It can be proved by induction on $i.$

{\rm (ii)} We show this by induction on $k$. It is clear for $k =1.$ Now assume $k >1$.  Then we have $\sigma (x_{i_1}, x_{i_2}, \cdots , x_{i_{2k}})= [\sigma _1(x_{i_1}, x_{i_2}, \cdots , x_{i_{2s}}), \sigma _2(x_{i_{2s+1}}, x_{i_2}, \cdots , x_{i_{2k}})]${\tiny $_{R}$} or $\sigma (x_{i_1}, x_{i_2}, \cdots , x_{i_{2k}})$
$ = [\sigma _1(x_{i_1}, x_{i_2}, \cdots , x_{i_{2s+1}}), \sigma _2(x_{i_{2s+2}}, x_{i_2}, \cdots , x_{i_{2k}})]$
{\tiny $_{R}$}. In the first case, we have

$\sigma (x_{i_1}, x_{i_2}, \cdots , x_{i_{2k}}) = [\sigma _1(x_{i_1}, x_{i_2}, \cdots , x_{i_{2s}}), \sigma _2(x_{i_{2s+1}}, x_{i_2}, \cdots , x_{i_{2k}})]${\tiny $_{R}$}

$= [ \alpha_s ( p_{1,2}^{s}(x_1x_2)^s- p_{2,1}^{s} (x_2x_1)^s), \alpha_t ( p_{1,2}^{t}(x_1x_2)^t-  p_{2,1}^{t}(x_2x_1)^t) ]
${\tiny $_{R}$}

\hfill (by inductive assumption, where $ s +t =k$)

$= \alpha_k (p_{1,2}^{k} (x_1x_2)^k- p_{2,1}^{k} (x_2x_1)^k)$(where $ \alpha_k=0$ ).

\noindent
In second case, we have

$\sigma (x_{i_1}, x_{i_2}, \cdots , x_{i_{2k}}) = [\sigma _1(x_{i_1}, x_{i_2}, \cdots , x_{i_{2s+1}}), \sigma _2(x_{i_{2s+2}}, x_{i_2}, \cdots , x_{i_{2k}})]$
{\tiny $_{R}$}

$= [ \beta _s ( x_{j_1}x_{j_2})^s x_{j_1}, \beta _t (x_{k_1}x_{k_2})^t x_{k_1})]$
{\tiny $_{R}$}

\hfill (by inductive assumption, where $ s + t+1 =k,  j_1 \not=j_2 $ and $ k_1\not= k_2 $)

$=\left \{  \begin{array}{ll}
0,  & \mbox {when }    k_{1}=j_{1} \\
\beta_{s} \beta _{t}(-1)^{s+t}(p_{j_{1},j_{2}}
p_{j_{2},j_{1}})^{st}( p_{j_{1},j_{2}}^{k}(x_{j_{1}}x_{j_{2}})^{k}-  p_{j_{2},j_{1}}^{k}(x_{j_{2}}x_{j_{1}})^{k}),   & \mbox {when }  k_{1}\neq j_{1}  \\
\end{array}\right. $

$= \alpha_k (p_{1,2}^{k} (x_1x_2)^k- p_{2,1}^{k} (x_2x_1)^k)$.

$\sigma (x_{i_1}, x_{i_2}, \cdots , x_{i_{2k+1}}) = b[\sigma _1(x_{i_1}, x_{i_2}, \cdots , x_{i_{2s}}), \sigma _2(x_{i_{2s+1}}, x_{i_{2s+2}}, \cdots , x_{i_{2k+1}})]${\tiny $_{R}$}

\hfill (where $b=1$ or $-1$)

$= [\alpha_s (p_{1,2}^{s} (x_1x_2)^s- p_{2,1}^{s} (x_2x_1)^s), \beta _t (x_{k_1}x_{k_2})^t x_{k_1})]${\tiny $_{R}$}

\hfill (by inductive assumption, where $s+t =k $ and $ k_1\not= k_2 $)

$=\left \{  \begin{array}{ll}
2\alpha_s \beta _t(-1)^{s}(p_{1,2}p_{2,1})^{st+s}(x_{1}x_{2})^{k}x_{1},  & \mbox {when }    k_{1}=1 \\
2\alpha_s \beta _t(-1)^{s+1}(p_{1,2}p_{2,1})^{st+s}(x_{2}x_{1})^{k}x_{2},   & \mbox {when }  k_{1}= 2  \\
\end{array}\right. $

$= \beta_k (x_{k_1}x_{k_2})^kx_{k_1}$.

Consequently, {\rm (ii)} holds.

{\rm (iii)}
We show this  by induction on $k$.

$
[[ [[x_{1}, x_{2}]${\tiny $_{R}$}$, x_1 ]${\tiny $_{R}$}$, x_2]${\tiny $_{R}$}$\cdots , x_{1}]${\tiny $_{R}$}$]${\tiny $_{R}$}

$= 2^{k-2}
 (p_{1,2}p_{2,1})^{\frac{(k-1)(k-2)}{2}} [[ (p_{1,2}^{k-1}(x_{1}x_{2})^{k-1}- p_{2,1}^{k-1}(x_{2}x_{1})^{k-1}   ), x_1]${\tiny $_{R}$}
$, x_2 ]${\tiny $_{R}$}

\hfill (by inductive assumption)

$= 2^{k-2}
 (p_{1,2}p_{2,1})^{\frac{(k-1)(k-2)}{2}} [2(-1)^{k-1} (p_{1,2}p_{2,1})^{k-1}(x_{1}x_{2})^{k-1}x_{1}, x_2 ]${\tiny $_{R}$}

$= (-1)^{k-1}2^{k-1}
 (p_{1,2}p_{2,1})^{\frac{k(k-1)}{2}} ((-1)^{k-1}p_{1,2}^{k}(x_{1}x_{2})^{k}
 -(-1)^{k-1}p_{2,1}^{k}(x_{2}x_{1})^{k})$

$=2^{k-1}
 (p_{1,2}p_{2,1})^{\frac{k(k-1)}{2}}(p_{1,2}^{k} (x_{1}x_{2})^{k}- p_{2,1}^{k}(x_{2}x_{1})^{k}   )$. \qed

\begin {Theorem}  \label {4.3}   Assume that $\mathfrak L(V)${\tiny $_{R}$} is a connected Nichols R-bicharacter algebra of diagonal type with $\dim V=2$ and $p_{1,1}=p_{2,2}= -1$. Let
$P_{1}:=$

\noindent $\left \{  \begin{array}{ll}
\{p_{1,2}^{k+1}(x_{1}x_2)^{k+1}-p_{2,1}^{k+1}(x_{2}x_1)^{k+1},x_2(x_{1}    x_2)^{k},(x_{1}  x_2)^{k}x_1,  0 \le k <  m \}\\
\hfill  -\{p_{1,2}^{m}(x_{1}x_2)^{m}-p_{2,1}^{m}(x_{2}x_1)^{m }\}\\
\hfill   \hbox {when }
   (-p_{1,2}^{-2}p_{2,1})^{m}= (-p_{1,2}p_{2,1}^{-2})^{m}=-1 \hbox { with } {\rm ord}
(p_{1,2} p_{2,1})  =m < \infty \\
\{p_{1,2}^{k+1}(x_{1}x_2)^{k+1}-p_{2,1}^{k+1}(x_{2}x_1)^{k+1},x_2(x_{1}    x_2)^{k},(x_{1}  x_2)^{k}x_1,  0 \le k <  {\rm ord}
(p_{1,2} p_{2,1})  \}, \hfill  \hbox {otherwise. } \\
\end{array}\right. $

Then $P_{1}$ is a basis of $\mathfrak L(V)${\tiny $_{R}$}. Furthermore, if  $ {\rm ord}
(p_{1,2} p_{2,1}) = m < \infty$, then

$\dim \mathfrak L(V)${\tiny $_{R}$} $=\left \{  \begin{array}{ll}
3m-1,  &    \hfill  \hbox {when }
   (-p_{1,2}^{-2}p_{2,1})^{m}= (-p_{1,2}p_{2,1}^{-2})^{m}=-1,  \\
3m,   & \mbox {otherwise. } \\
\end{array}\right. .$

\end {Theorem}
\noindent {\it Proof.} By Lemma \ref {4.1}(ii), Lemma \ref {4.2} (iv), $P_{1}$ is linearly independent.
It follows from  Lemma \ref {4.2} (i) (iii) that $P_{1} \subseteq   \mathfrak L(V)${\tiny $_{R}$}.
By Lemma \ref {4.2} (ii),  $P_{1}$ is a basis of $\mathfrak L(V)${\tiny $_{R}$}. \qed

\subsection{Nichols L-bicharacter algebra}

\begin {Lemma}\label {4.4}   Assume that $\mathfrak L(V)${\tiny $_{L}$} is a connected Nichols L-bicharacter algebra of diagonal type with $\dim V=2$ and $p_{1,1}=p_{2,2}= -1$.

{\rm (i)} $(r_{x_{1}}r_{x_{2}})^{i}[x_{1}]${\tiny $_{L}$}$ = (-2)^{i} (p_{1,2}p_{2,1})^{\frac{i(i+1)}{2}}(x_{1}x_{2})^{i}x_{1}$ and

$(r_{x_{2}}r_{x_{1}})^{i}[x_{2}]${\tiny $_{L}$}$ = (-2)^{i} (p_{1,2}p_{2,1})^{\frac{i(i+1)}{2}}x_{2}(x_{1}x_{2})^{i}$ for $i \ge 0.$

{\rm (ii)} For any method $\sigma$ adding bracket  $[\cdot,\cdot]${\tiny $_{L}$} and  $k >0$, there exist $\alpha_k, \beta _k \in F$
such that
$\sigma (x_{i_1}, x_{i_2}, \cdots , x_{i_{2k+1}}) = \beta_k  (x_1x_2)^ {k} x_1$ or $ \beta_k  x_2(x_1x_2)^ {k}$
and $\sigma (x_{i_1}, x_{i_2}, \cdots , x_{i_{2k}}) = \alpha_k  ( p_{2,1}^{k}(x_1x_2)^ k - p_{1,2}^{k}(x_2x_1)^ k)$.

{\rm (iii)} $
[[ [[x_{1}, x_{2}]${\tiny $_{L}$}$, x_1 ]${\tiny $_{L}$}$, x_2]${\tiny $_{L}$}$\cdots , x_{1}]${\tiny $_{L}$}$]${\tiny $_{L}$}$ = 2^{k-1}
 (p_{1,2}p_{2,1})^{\frac{k(k-1)}{2}}(p_{2,1}^{k} (x_{1}x_{2})^{k}- p_{1,2}^{k}(x_{2}x_{1})^{k}   )$,
where there exist $2k-1$ brackets in the left hand side.

{\rm (iv)} $$<y_1, p_{2,1}^{i}(x_{1}x_2)^{i}-p_{1,2}^{i}(x_{2}x_1)^{i}>=(-1)^{i}
(1+(-p_{2,1})^{i})x_2(x_{1}    x_2)^{i-1},$$
$$<y_2, p_{2,1}^{i}(x_{1}x_2)^{i}-p_{1,2}^{i}(x_{2}x_1)^{i}>=(-1)^{i+1}
(1+(-p_{1,2})^{i})(x_{1}    x_2)^{i-1}x_1.$$
\end {Lemma}
\noindent {\it Proof.}
{\rm (i)} It can be proved by induction on $i.$

{\rm (ii)} We show this by induction on $k$. It is clear for $k =1.$ Now assume $k >1$.  Then we have $\sigma (x_{i_1}, x_{i_2}, \cdots , x_{i_{2k}})= [\sigma _1(x_{i_1}, x_{i_2}, \cdots , x_{i_{2s}}),$ $\sigma _2(x_{i_1}, x_{i_2}, \cdots , x_{i_{2t}})]${\tiny $_{L}$} or $\sigma (x_{i_1}, x_{i_2}, \cdots , x_{i_{2k}})$
$= [\sigma _1(x_{i_1}, x_{i_2}, \cdots , x_{i_{2s-1}}),$ $\sigma _2(x_{i_1}, x_{i_2}, \cdots , x_{i_{2t-1}})]${\tiny $_{L}$}. In the first case, we have

$\sigma (x_{i_1}, x_{i_2}, \cdots , x_{i_{2k}}) = [\sigma _1(x_{i_1}, x_{i_2}, \cdots , x_{i_{2s}}), \sigma _2(x_{i_{2s+1}}, x_{i_2}, \cdots , x_{i_{2k}})]${\tiny $_{L}$}

$= [ \alpha_s ( p_{2,1}^{s}(x_1x_2)^s- p_{12}^{s} (x_2x_1)^s), \alpha_t ( p_{2,1}^{t}(x_1x_2)^t-  p_{1,2}^{t}(x_2x_1)^t) ]
${\tiny $_{L}$}

\hfill (by inductive assumption, where $ s +t =k$)

$= \alpha_k (p_{2,1}^{k} (x_1x_2)^k- p_{1,2}^{k} (x_2x_1)^k)$(where $ \alpha_k=0$ ).

\noindent
In second case, we have

$\sigma (x_{i_1}, x_{i_2}, \cdots , x_{i_{2k}}) = [\sigma _1(x_{i_1}, x_{i_2}, \cdots , x_{i_{2s+1}}), \sigma _2(x_{i_{2s+2}}, x_{i_2}, \cdots , x_{i_{2k}})]$
{\tiny $_{L}$}

$= [ \beta _s ( x_{j_1}x_{j_2})^s x_{j_1}, \beta _t (x_{k_1}x_{k_2})^t x_{k_1})]$
{\tiny $_{L}$}

\hfill (by inductive assumption, where $ s + t+1 =k,  j_1 \not=j_2 $ and $ k_1\not= k_2 $)

$=\left \{  \begin{array}{ll}
0,  & \mbox {when }    k_{1}=j_{1} \\
\beta_{s} \beta _{t}(-1)^{s+t}(p_{j_{1},j_{2}}
p_{j_{2},j_{1}})^{st}( p_{j_{2},j_{1}}^{k}(x_{j_{1}}x_{j_{2}})^{k}-  p_{j_{1},j_{2}}^{k}(x_{j_{2}}x_{j_{1}})^{k}),   & \mbox {when }  k_{1}\neq j_{1}  \\
\end{array}\right. $

$= \alpha_k (p_{2,1}^{k} (x_1x_2)^k- p_{1,2}^{k} (x_2x_1)^k)$.

$\sigma (x_{i_1}, x_{i_2}, \cdots , x_{i_{2k+1}}) = b[\sigma _1(x_{i_1}, x_{i_2}, \cdots , x_{i_{2s}}), \sigma _2(x_{i_{2s+1}}, x_{i_{2s+2}}, \cdots , x_{i_{2k+1}})]${\tiny $_{L}$}

\hfill (where $b=1$ or $-1$)

$= [\alpha_s (p_{2,1}^{s} (x_1x_2)^s- p_{1,2}^{s} (x_2x_1)^s), \beta _t (x_{k_1}x_{k_2})^t x_{k_1})]${\tiny $_{L}$}

\hfill (by inductive assumption, where $s+t =k $ and $ k_1\not= k_2 $)

$=\left \{  \begin{array}{ll}
2\alpha_s \beta _t(-1)^{s}(p_{1,2}p_{2,1})^{st+s}(x_{1}x_{2})^{k}x_{1},  & \mbox {when }    k_{1}=1 \\
2\alpha_s \beta _t(-1)^{s+1}(p_{1,2}p_{2,1})^{st+s}(x_{2}x_{1})^{k}x_{2},   & \mbox {when }  k_{1}= 2  \\
\end{array}\right. $

$= \beta_k (x_{k_1}x_{k_2})^kx_{k_1}$.

Consequently, {\rm (ii)} holds.

{\rm (iii)}
We show this  by induction on $k$.

$[[ [[x_{1}, x_{2}]${\tiny $_{L}$}$, x_1 ]${\tiny $_{L}$}$, x_2]${\tiny $_{L}$}$\cdots , x_{1}]${\tiny $_{L}$}$]${\tiny $_{L}$}

$= 2^{k-2}
 (p_{1,2}p_{2,1})^{\frac{(k-1)(k-2)}{2}} [[ (p_{2,1}^{k-1}(x_{1}x_{2})^{k-1}- p_{1,2}^{k-1}(x_{2}x_{1})^{k-1}   ), x_1]${\tiny $_{L}$}
$, x_2 ]${\tiny $_{L}$}

\hfill (by inductive assumption)

$= 2^{k-2}
 (p_{1,2}p_{2,1})^{\frac{(k-1)(k-2)}{2}} [2(-1)^{k-1} (p_{1,2}p_{2,1})^{k-1}(x_{1}x_{2})^{k-1}x_{1}, x_2 ]${\tiny $_{L}$}

$= (-1)^{k-1}2^{k-1}
 (p_{1,2}p_{2,1})^{\frac{k(k-1)}{2}} ((-1)^{k-1}p_{2,1}^{k}(x_{1}x_{2})^{k}
 -(-1)^{k-1}p_{1,2}^{k}(x_{2}x_{1})^{k} )$

$=2^{k-1}
 (p_{1,2}p_{2,1})^{\frac{k(k-1)}{2}}(p_{2,1}^{k} (x_{1}x_{2})^{k}- p_{1,2}^{k}(x_{2}x_{1})^{k}   )$. \qed

\begin {Theorem}  \label {4.5}   Assume that $\mathfrak L(V)${\tiny $_{L}$} is a connected Nichols L-bicharacter algebra of diagonal type with $\dim V=2$ and $p_{1,1}=p_{2,2}= -1$. Let
$P_{2}:=$

\noindent $\left \{  \begin{array}{ll}
\{p_{2,1}^{k+1}(x_{1}x_2)^{k+1}-p_{1,2}^{k+1}(x_{2}x_1)^{k+1},x_2(x_{1}    x_2)^{k},(x_{1}  x_2)^{k}x_1,  0 \le k <  m \}\\
\hfill  -\{p_{2,1}^{m}(x_{1}x_2)^{m}-p_{1,2}^{m}(x_{2}x_1)^{m }\}\\
\hfill   \hbox {when }
   (-p_{2,1})^{m}= (-p_{1,2})^{m}=-1 \hbox { with } {\rm ord}
(p_{1,2} p_{2,1})  =m < \infty \\
\{p_{2,1}^{k+1}(x_{1}x_2)^{k+1}-p_{1,2}^{k+1}(x_{2}x_1)^{k+1},x_2(x_{1}    x_2)^{k},(x_{1}  x_2)^{k}x_1,  0 \le k <  {\rm ord}
(p_{1,2} p_{2,1})  \}, \hfill  \hbox {otherwise. } \\
\end{array}\right. $

Then $P_{2}$ is a basis of $\mathfrak L(V)${\tiny $_{L}$}. Furthermore, if  $ {\rm ord}
(p_{1,2} p_{2,1}) = m < \infty$, then

\noindent $\dim \mathfrak L(V)${\tiny $_{L}$} $=\left \{  \begin{array}{ll}
3m-1,  &    \hfill  \hbox {when }
   (-p_{2,1})^{m}= (-p_{1,2})^{m}=-1,  \\
3m,   & \mbox {otherwise. } \\
\end{array}\right. .$

\end {Theorem}
\noindent {\it Proof.} By Lemma \ref {4.1}(ii),Lemma \ref {4.4} (iv), $P_{2}$ is linearly independent.
It follows from  Lemma \ref {4.4} (i) (iii) that $P_{2} \subseteq   \mathfrak L(V)${\tiny $_{L}$}.
By Lemma \ref {4.4} (ii),  $P_{2}$ is a basis of $\mathfrak L(V)${\tiny $_{L}$}. \qed

\section {Some equivalent characterizations}\label {s4}
In this section we give the sufficient and necessary conditions for $\mathfrak L(V)${\tiny $_{R}$}$= \mathfrak L(V)$, $\mathfrak L(V)${\tiny $_{L}$}$= \mathfrak L(V)$, $\mathfrak B(V) = F\oplus \mathfrak L(V)${\tiny $_{R}$} and $\mathfrak B(V) = F\oplus \mathfrak L(V)${\tiny $_{L}$}, respectively.

\begin{Lemma}\label{3.3} Assume that   $u_{i}$ is a homogeneous element and  $u_{i}u_{j} = p_{u_{i}, u_{j}} u_{j}u_{i}$  with $p_{u_{i}, u_{j}} \in F^*$ and $p_{u_{i}, u_{j}}p_{u_{i}, u_{j}}=1$ for $1\le i ,  j \le k$ and   $p_{u_{i}, u_{j}} =1$ when $u_{i}= u_{j}$.  Then

$[u_{1 }, \cdots,  u_{m }]${\tiny $_{R}$}$= \prod \limits _{j=1}^{m-1}p_{u_{m }\cdots u_{j+1 },u_{j }}( p_{u_{j },u_{m }\cdots u_{j+1 }}^{3}-1)u_{m }\cdots u_{2 } u_{1}$,

\noindent where $[u_{1 }, \cdots,  u_{m }]${\tiny $_{R}$}$:= [u_1, [ u_2, \cdots [u_{m-1}, u_m]${\tiny $_{R}$}$ \cdots ]${\tiny $_{R}$}$]${\tiny $_{R}$}.
\end {Lemma}
\noindent {\it Proof.} We show  Lemma by induction on $m$. Obviously, $u_{1 }[u_{2 }, \cdots,  u_{m }]${\tiny $_{R}$}  $=p_{u_{1 }, u_{2 }\cdots u_{m }}[u_{2 }\cdots u_{m }]${\tiny $_{R}$} $u_{1 }$. Therefore,
$[u_{1 }, \cdots,  u_{m }]${\tiny $_{R}$} $=p_{u_{1 }, u_{2 }\cdots u_{m }}u_{1 } [u_{2 }, \cdots,  u_{m }]${\tiny $_{R}$}$-p_{u_{2 }\cdots u_{m },u_{1 }}[u_{2 }, \cdots,  u_{m }]${\tiny $_{R}$}$u_{1 }$

$=(p_{u_{1 }, u_{2 }\cdots u_{m }}^{2}-p_{u_{2 }\cdots u_{m },u_{1 }} )[u_{2 }, \cdots,  u_{m }]${\tiny $_{R}$}$u_{1 }$

$= \prod \limits _{j=1}^{m-1}p_{u_{m }\cdots u_{j+1 },u_{j }}(p_{u_{j },u_{m }\cdots u_{j+1 }}^{3}-1)u_{m }\cdots u_{2 } u_{1}$ (by inductive assumption). \qed

\begin {Lemma} \label {3.4}  Let $(L_{0} , [\cdot,\cdot]${\tiny $_{R}$}$)$ be a Lie algebra and $u_1, u_2, \cdots, u_m \in L_{0}$.
If $\sigma$ is a  method  of adding bracket $[\cdot,\cdot]${\tiny $_{R}$} on $u_1,  u_2,  \cdots,  u_m$,  then there exist some $\tau_{j}\in \mathbb S_{m}$, $\xi _{j}=1$ or $-1$ such that
$$(\ast)~~~~\sigma (u_1, u_2, \cdots, u_m)=\sum\limits_{j=1} ^r \xi _{j}[u_{\tau_{j}(1)}, \cdots,  u_{\tau_{j}(m)}]_{\tiny R}.$$

\end {Lemma}
 \noindent {\it Proof.}  We show  $(\ast)$ by induction on $m$. Obviously,  $(\ast)$ holds for $m =2$. Assume $m>2.$
 Let $\sigma (u_1, \cdots, u_m) = [\sigma_1 (u_1,  \cdots, u_s), $ $\sigma _2(u_{s+1}, u_{s+2}, $ $\cdots,$ $ u_m)]${\tiny $_{R}$}.

Now we show $(\ast)$ by induction on $s$. In case $s=1$,

$\sigma (u_1, \cdots, u_m) = [u_1, \sigma _2(u_{2}, \cdots, u_m)]
${\tiny $_{R}$}

$=
 \sum\limits_{j=1} ^r \xi _{j} [ u_1, [u_{\tau_{j}(2)},\cdots, u_{\tau_{j}(m)}]
 ${\tiny $_{R}$}$ ]${\tiny $_{R}$}
(where $\tau_j \in \mathbb S_{ \{2, 3. \cdots, m \}}$ for  $1\le j \le r) $

$= \sum\limits_{j=1} ^r \xi _{j} [ u_{\tau _j (1)}, u_{\tau_{j}(2)},\cdots, u_{\tau_{j}(m)}]${\tiny $_{R}$}.  Therefore, (*) holds.

Assume $s>1$ and
$\sigma_1 (u_1, \cdots, u_s) = [\sigma_3 (u_1, \cdots, u_k)$, $\sigma _4 (u_{k+1}, $ $u_{k+2}, $ $\cdots, u_s)]${\tiny $_{R}$}.

See  $\sigma (u_1, \cdots, u_m) = [ [\sigma_3 (u_1, \cdots, u_k), \sigma _4 (u_{k+1}, u_{k+2}, \cdots, u_s)]${\tiny $_{R}$}$  , \sigma _2(u_{s+1}, u_{s+2}, \cdots, u_m)]${\tiny $_{R}$}

$=[ [\sigma_3 (u_1, \cdots, u_k), \sigma _2 (u_{s+1}, u_{s+2}, \cdots, u_m)]${\tiny $_{R}$}$ , \sigma _4 (u_{k+1}, u_{k+2}, \cdots, u_s)]${\tiny $_{R}$}

\hfill $+[ \sigma_3 (u_1, \cdots, u_k), [ \sigma _4 (u_{k+1}, u_{k+2}, \cdots, u_s),  \sigma _2 (u_{s+1}, u_{s+2}, \cdots, u_m)]${\tiny $_{R}$}$ ]${\tiny $_{R}$}

$=\sum\limits_{j=1} ^{r_1} \xi _{j} [[ u_ {\tau _j (1)}, \cdots, u_{\tau_{j}(k)}, u_ {\tau _j (s+1)}, u_{\tau_{j}(s+ 2)},\cdots, u_{\tau_{j}(m)}]${\tiny $_{R}$}$, \sigma _4 (u_{k+1}, u_{k+2}, \cdots, u_s)]${\tiny $_{R}$}

$+
 \sum\limits_{j=r_1+1} ^{r_2} [ \sigma_3 (u_1, \cdots, u_k),  \xi _{j} [ u_ {\tau _j (k+1)}, u_{\tau_{j}(k+2)},\cdots, u_{\tau_{j}(s)}, u_ {\tau _j (s+1)}, u_{\tau_{j}(s+ 2)},\cdots, u_{\tau_{j}(m)}]${\tiny $_{R}$}$]${\tiny $_{R}$}

(by first inductive assumption, where $ \tau_j \in \mathbb S_{ \{1,2, 3. \cdots, k, s+1, s+2, \cdots, m \}}$  for $ 1\le j \le r_1$;
 $\tau_j \in \mathbb S_{ \{k +1, k+2. \cdots,  m \}}$ for $ r_1 +1\le j \le r_2)$

$=\sum\limits_{j=r_2+ 1} ^{r_3} \xi _{j}[u_{\tau_{j}(1)}, \cdots,  u_{\tau_{j}(m)}]${\tiny $_{R}$}

(by second inductive assumption, where  $\tau_j \in \mathbb S_{ \{1, 2. \cdots,  m \}}$ for $ r_2+1\le j \le r_3$).

Therefore, $(\ast)$ holds. \qed

\begin {Proposition} \label {5.1} If $\mathfrak B(V) $ is a Nichols algebra of diagonal type, then $\mathfrak L(V)= \mathfrak L(V)${\tiny $_{R}$} if and only if $p_{i,i}^2=1,p_{i,j}p_{j,i}=1,p_{i,j}\in R_{3}$ for $ 1\leq i\neq j\leq n$. In this case, $\mathfrak L(V)= \mathfrak L(V)${\tiny $_{R}$}$=V$.
\end {Proposition}
\noindent {\it Proof.} The sufficiency. By \cite [Coroll.2.5] {WWZZ18}, $\mathfrak L(V) = V$. By Lemma \ref {1.4}, $[x_i, x_j]${\tiny $_{R}$}$ =0$ for $i \not= j$ and  $\mathfrak L(V)${\tiny $_{R}$}$ = V$.

The necessity. We show this by following three steps.

{\rm (i)} $p_{i,i}^2 =1$ since $x_i ^k \in \mathfrak L(V)$ for $k \le {\rm ord} (p_{i,i})$ and $x_i ^m \notin \mathfrak L(V)${\tiny $_{R}$} when $0\not=x_i ^m$ and $m>1.$

{\rm (ii)}   If $p_{i,j} p_{j,i} =1$ and $p_{i,j} \notin R_{3}$ with $i \not= j,$ then $[x_i, x_j]${\tiny $_{R}$}$ \not=0$ and $[x_i, x_j]=0$. Consequently, $[x_i, x_j]${\tiny $_{R}$}$ \in \mathfrak L(V)${\tiny $_{R}$}$ - \mathfrak L(V)$, which is a contradiction.

{\rm (iii)} If $p_{i,j} p_{j,i} \not=1$ with $i < j$, then $0\not= [x_i, x_j] \in \mathfrak L(V)= \mathfrak L(V)${\tiny $_{R}$} and $[x_i, x_j] = k[x_i, x_j] ${\tiny $_{R}$} with $k \in F^{*}$. By
\cite [Coroll.2.5] {WWZZ18}, $x_jx_i \in \mathfrak L(V)$, which implies $x_jx_i = k'[x_i, x_j]${\tiny $_{R}$} with $k' \in F^{*}$. This is a contradiction since $x_jx_i$ and $[x_i, x_j]$ are linearly  independent. \qed

\begin {Proposition} \label {5.2} If $\mathfrak B(V) $ is a Nichols algebra of diagonal type, then the following conditions are equivalent:

(1) $\mathfrak L(V)= \mathfrak L(V)${\tiny $_{L}$}.
(2) $\mathfrak L(V)= \mathfrak L^{-}(V)$.

(3) $p_{i,i}^2=1,p_{i,j}=p_{j,i}=1$ for $ 1\leq i\neq j\leq n$. In this case, $\mathfrak L(V)= \mathfrak L(V)${\tiny $_{L}$}$=V$.
\end {Proposition}
\noindent {\it Proof.} By \cite [Prop.6.3] {WWZZ18}, $(2)\Longleftrightarrow(3)$. The proof of $(1)\Longleftrightarrow(3)$ is similar to the proof of $(2)\Longleftrightarrow(3)$. \qed

\begin {Proposition} \label {5.3} Assume that $\mathfrak B(V) $ is a Nichols algebra of diagonal type. Then $\mathfrak B(V) = F \oplus \mathfrak L(V)${\tiny $_{R}$} if and only if $p_{i,i}=-1,p_{i,j}p_{j,i}=1$ for all $1\leq i\neq j\leq n$ and  there  exist $\tau\in \mathbb  S_{m}$ such that $\prod \limits _{j=1}^{m-1}
(p_{h_{\tau(j)},h_{\tau(m)}\cdots h_{\tau(j+1)}}^{3}-1)\neq0$ for
all $h_{1}>h_{2}>\cdots>h_{m}$ with $h_{i}\in \{x_{1},\ldots,x_{n}\},1\leq i\leq m$.
\end {Proposition}
\noindent {\it Proof.} The necessity. If there exists $1\le i \le n$ such that  $ p_{i,i}\neq-1$, then $0\neq x_{i}^{2}\in \mathfrak B(V)$ and $x_{i}^{2}\notin \mathfrak L(V)${\tiny $_{R}$}, which  is a contradiction. If there exist $i,  j $ such that $ p_{i,j}p_{j,i}\neq1$ with  $1\le i< j \le n$, then $[x_{i},x_{j}]\neq0$ and $[x_{i},x_{j}]${\tiny $_{R}$}$\neq0$. Since $\mathfrak B(V) = F \oplus \mathfrak L(V)${\tiny $_{R}$}, we have that there exist  $ k, k'\in F^{*}$ such that $[x_{i},x_{j}]=k[x_{i},x_{j}]${\tiny $_{R}$}, and $x_{j}x_{i}=k'[x_{i},x_{j}]${\tiny $_{R}$}, which contradicts  to that $[x_{i},x_{j}]$ and $x_{j}x_{i}$ are linearly independent. Therefore, $V$ is a quantum linear space.

If there exist $ \ h_{1}>h_{2}>\cdots>h_{m}$ with $h_{i}\in \{x_{1},\ldots,x_{n}\},1\leq i\leq m$,  such that $\prod \limits _{j=1}^{m-1}( p_{h_{\tau(j)},h_{\tau(m)}\cdots h_{\tau(j+1)}}^{3}-1)= 0$ for any $\tau \in \mathbb S_m$.  By \cite [Lem.3.2] {WWZZ18}, Lemma \ref {3.3} and Lemma \ref {3.4}, $\sigma ( h_{\tau (1)}, h_{\tau (2)}, \cdots, h_{\tau (m)} ) =0$ for any $\tau \in \mathbb S_m$ and  any  method $\sigma$ of adding bracket $[\cdot,\cdot]${\tiny $_{R}$} on $ h_{\tau (1)}, h_{\tau (2)}, \cdots, h_{\tau (m)}$. Consequently, $0\not= h_1h_2\cdots h_m \notin \mathfrak L(V)${\tiny $_{R}$}, which is a contradiction.

The sufficiency. Obviously, $V$ is a quantum linear space. For $ \forall\ h_{1}>h_{2}>\cdots>h_{m}$ with $h_{i}\in \{x_{1},\ldots,x_{n}\},1\leq i\leq m$,
there  exist $\tau\in \mathbb  S_{m}$ such that $\prod \limits _{j=1}^{m-1}( p_{h_{\tau(j)},h_{\tau(m)}\cdots h_{\tau(j+1)}}^{3}-1)\neq0$. By Lemma \ref {3.3}, we have that there exist  $a\in F^{*}$ such that $h_1h_2\cdots h_m =a h _{\tau (m)} h _{\tau (m-1)} $ $ \cdots $ $h _{\tau (1)} \in \mathfrak  L(V)${\tiny $_{R}$}. \qed

\begin{Lemma}\label{3.8} Assume that   $u_{i}$ is a homogeneous element and  $u_{i}u_{j} = p_{u_{i}, u_{j}} u_{j}u_{i}$  with $p_{u_{i}, u_{j}} \in F^*$ and $p_{u_{i}, u_{j}}p_{u_{i}, u_{j}}=1$ for $1\le i ,  j \le k$ and   $p_{u_{i}, u_{j}} =1$ when $u_{i}= u_{j}$.  Then

$[u_{1 }, \cdots,  u_{m }]${\tiny $_{L}$}$= \prod \limits _{j=1}^{m-1}p_{u_{m }\cdots u_{j+1 },u_{j }}( p_{u_{j },u_{m }\cdots u_{j+1 }}-1)u_{m }\cdots u_{2 } u_{1}$,

\noindent where $[u_{1 }, \cdots,  u_{m }]${\tiny $_{L}$}$:= [u_1, [ u_2, \cdots [u_{m-1}, u_m]${\tiny $_{L}$}$ \cdots ]${\tiny $_{L}$}$]${\tiny $_{L}$}.
\end {Lemma}

\begin {Lemma} \label {3.9}  Let $(L_{0} , [\cdot,\cdot]${\tiny $_{L}$}$)$ be a Lie algebra and $u_1, u_2, \cdots, u_m \in L_{0}$.
If $\sigma$ is a  method  of adding bracket $[\cdot,\cdot]${\tiny $_{L}$} on $u_1,  u_2,  \cdots,  u_m$,  then there exist some $\tau_{j}\in \mathbb S_{m}$, $\xi _{j}=1$ or $-1$ such that
$(\ast)\ \ \ \sigma (u_1, u_2, \cdots, u_m)=\sum\limits_{j=1} ^r \xi _{j}[u_{\tau_{j}(1)}, \cdots,  u_{\tau_{j}(m)}]${\tiny $_{L}$}.
\end {Lemma}

\begin {Proposition} \label {5.4} Assume that $\mathfrak B(V) $ is a Nichols algebra of diagonal type. Then the following conditions are equivalent:

(1) $\mathfrak B(V) = F \oplus \mathfrak L(V)${\tiny $_{L}$}.
(2) $\mathfrak B(V) = F \oplus \mathfrak L^{-}(V)$.

(3) $p_{i,i}=-1,p_{i,j}p_{j,i}=1$ for all $1\leq i\neq j\leq n$ and  there  exist $\tau\in \mathbb  S_{m}$ such that $\prod \limits _{j=1}^{m-1}
(p_{h_{\tau(j)},h_{\tau(m)}\cdots h_{\tau(j+1)}}-1)\neq0$ for
all $h_{1}>h_{2}>\cdots>h_{m}$ with $h_{i}\in \{x_{1},\ldots,x_{n}\},1\leq i\leq m$.
\end {Proposition}
\noindent {\it Proof.} By \cite [Prop.6.4] {WWZZ18}, $(2)\Longleftrightarrow(3)$. The proof of $(1)\Longleftrightarrow(3)$ is similar to the proof of $(2)\Longleftrightarrow(3)$. \qed

\begin {Conjecture} \label {5.8} Assume that $\mathfrak B(V) $ is a Nichols algebra of diagonal type. Then $\mathfrak L(V)${\tiny $_{L}$}$=\mathfrak L^{-}(V)$.
\end {Conjecture}

\begin {Question} \label {5.9} Assume that $\mathfrak B(V) $ is a Nichols algebra of diagonal type. Give the sufficient and necessary conditions for $\mathfrak L(V)${\tiny $_{L}$}$=\mathfrak L(V)${\tiny $_{R}$}.
\end {Question}

\section {Classification of $\mathfrak L(V)${\tiny $_{L}$} and $\mathfrak L(V)${\tiny $_{R}$}}\label {s5}
In this section it is proved that if $\mathfrak B(V)$ is a connected Nichols algebra of diagonal type with $\dim V>1$, then $\mathfrak B(V)$ is finite-dimensional if and only if $\mathfrak L(V)${\tiny $_{L}$} is finite-dimensional if and only if $\mathfrak L(V)${\tiny $_{R}$} is finite-dimensional.

Let $| u| $ denote the length of word $u$.
A word $u\in W$ is called a Lyndon word if $| u| =1$ or $| u| \geq2$,  and for each
representation $u=u_1u_2$,  where $u_1$and $u_2$ are nonempty
words,  the inequality $u<u_2u_1$ holds (see \cite [Def. 1] {Kh99}).
Any word $u\in W$ has a unique decomposition into the product of
non-increasing sequence of Lyndon words by \cite [Th.5.1.5] {Lo83}.
If $u$ is a Lyndon word with $| u | >1$,  then there uniquely exist two Lyndon words $v$ and $w$
such that $u =vw$ and $v$ is shortest (see \cite [Prop. 5.1.3]{Lo83})(the composition is called the Shirshov decomposition of $u$).

We call $u$ is a standard   word with respect to $\mathfrak B(V)$ if $u $ can not be  written as a linear combination of  strictly greater  words in $\mathfrak B(V)$.

Let {\rm S}$ (\mathfrak B(V)):=\{ u \in W \mid  u \hbox { is a standard word with respect to } \mathfrak B(V) \}$, written as {\rm S} in short;
{\rm L } $:= \{ u\in W \mid u \hbox { is a Lyndon word}\}$.
Let {\rm H} $: =\{ u  \in {\rm L}\mid [u] \hbox { is a hard super-letter}\}$.
Notice that we view {\rm L} and {\rm H} are in $\mathfrak B(V)$ often for convenience.

Let $D =: \{[u] \mid [u] \hbox { is a hard super-letter}\}$.
$\Delta ^+(\mathfrak B(V)): =  \{ \deg (u) \mid [u]\in D\}$.
$ \Delta (\mathfrak B(V)) := \Delta ^+(\mathfrak B(V)) \cup \Delta ^-(\mathfrak B(V))$, which is called the root system of $V.$ If $ \Delta (\mathfrak B(V))$ is finite, then it is called an arithmetic root system.

\subsection {$\Delta (\mathfrak B(V))$ is not  an arithmetic root system}
\begin {Lemma}\label {6.1}{\rm (See \cite [Lem.2.2]{WZZ})}
{\rm (i)}  {\rm S} is a basis of $\mathfrak B(V)$.

{\rm (ii)} Any factor of a standard  word is a standard  word.

{\rm (iii)} If $u$ is a standard word, then  $u = u_{1}u_2 \cdots u_r$ with $u_1 \ge u_2\ge \cdots \ge u_r $ and $u_i \in {\rm S } \cap {\rm L}$ for $1\le i \le r.$
\end {Lemma}

\begin {Lemma}\label {6.2}{\rm (See \cite [Theo.2.4]{WZZ}) }  If $\mathfrak B(V)$ is a Nichols algebra of diagonal type, then ${\rm S} (\mathfrak B(V))\cap {\rm L} = {\rm H} (\mathfrak B(V))$.
\end {Lemma}

\begin {Lemma}\label {6.3}
{\rm (i)} If $l\in {\rm L}$, then $ [l]${\tiny $_{L}$}$=a_{l}  l + \sum \limits _ {w > l, \mid l \mid = \mid w\mid} a_{w} w$ in $\mathfrak B(V)$, where $a_l, a_w \in F$ with $a_l\not=0.$

{\rm (ii)} If $l \in  {\rm L}, $ then $ [l]${\tiny $_{R}$}$=b_{l} l + \sum \limits_{w > l, \mid l \mid = \mid w\mid} b_{w} w$ in $\mathfrak B(V)$, where $b_{l}, b_{w} \in F$ with $b_l\not=0.$
\end {Lemma}
\noindent {\it Proof.}  {\rm (i)} We show this by induction on $\mid l \mid.$  It is clear when $\mid l \mid =1$ since $[l]${\tiny $_{L}$}$=l$. Assume that $l = uv$ is the Shirshov decomposition of $l$.
If $u' > u$ and  $v'>v$ with $\mid u'\mid = \mid u\mid $
and $\mid v'\mid = \mid v\mid $, then $u'v' > uv =l$ and $v'u' >vu >l.$
\begin {eqnarray*} [l]{\tiny _{L}}&=&p _{v, u}[u][v] -p _{u, v} [v] [u]\\
 &=&p _{vu}(a_u 'u + \sum \limits _ {u' > u, \mid u' \mid = \mid u \mid} a_{u'}' u')
 (a_v' v + \sum \limits _ {v' > v, \mid v' \mid = \mid v \mid} a_{v'}' v') \\
&&- p _{u, v} (a_v' v + \sum \limits _ {v' > v, \mid v' \mid = \mid v \mid} a_{v'}' v')(a_u 'u + \sum \limits _ {u' > u, \mid u' \mid = \mid u \mid} a_{u'}' u')\ \ ( \hbox {by inductive  hypothesis}) \\
&=&  a_ll + \sum \limits _ {w > l, \mid l \mid = \mid w\mid} a_ww.
\end {eqnarray*}

{\rm (ii)}  The proof is similar to the proof of {\rm (i)}. \qed

\begin {Theorem}\label {6.4} If $\mathfrak B(V)$ is a Nichols algebra of diagonal type and $\Delta (\mathfrak B(V))$ is not  an arithmetic root system, then $\dim \mathfrak L(V)${\tiny $_{L}$}$ = \infty$,$\dim \mathfrak L(V)${\tiny $_{R}$}$ = \infty$.
\end {Theorem}

\subsection {$\Delta (\mathfrak B(V))$ is an arithmetic root system}
For our study of Nichols algebras we will need some non-standard formulas for quantum integers and Gaussian binomial coefficients.

In the ring $\mathbb Z[a]$, let $(0)_{a}=0$ and for any $m\in \mathbb N$,
$(m)_{a}=1+a+a^{2} +\cdots+a^{m-1}$. The polynomials $(m)_{a}$ with $m\in \mathbb Z$ are also
known as quantum integers. Moreover, let $(0)_{a}^{!}=1$, and for any $m\in \mathbb Z$ let $(m)_{a}^{!}=\prod \limits _{k=1}^{m}(k)_{a}$. For any $k,m\in \mathbb Z$ with $0\leq k\leq m$, the rational function
${\scriptsize\left(\begin{array}{cc} m\\ k \end{array}\right)_{a}}=\frac{(m)_{a}^{!}}{(k)_{a}^{!}(m-k)_{a}^{!}},$
is in fact an element of $\mathbb Z[a]$ and is called a Gaussian binomial coefficient.
For $m\in\mathbb N_0$,$k\in \mathbb Z$ with $k<0$ or $k>m$ one defines
${\scriptsize\left(\begin{array}{cc} m\\
k \end{array}\right)_{a}}=0$. The
Gaussian binomial coefficients satisfy the following formulas:
\begin {eqnarray}\label {e1}{\scriptsize\left(\begin{array}{cc} m\\
k \end{array}\right)_{a}}={\scriptsize\left(\begin{array}{cc} m\\
m-k \end{array}\right)_{a}},\end {eqnarray}
\begin {eqnarray}\label {e2}a^{k-1}{\scriptsize\left(\begin{array}{cc} m\\
k \end{array}\right)_{a}}+a^{m}{\scriptsize\left(\begin{array}{cc} m\\
k-1 \end{array}\right)_{a}}=a^{k-1}{\scriptsize\left(\begin{array}{cc} m+1\\
k \end{array}\right)_{a}}\end {eqnarray}
\noindent for $m\in \mathbb N,1\leq k\leq m$.

\begin {Lemma} \label {6.5} For $m\in \mathbb N,i\neq j$.

{\rm (i)} $\sum \limits _{k=0}^{m} ( - 1)^{k} {\scriptsize\left(       \begin{array}{c} m\\
k\end{array}\right)}(m-k)_{a}=(a-1)^{m-1}$.

{\rm (ii)} $l_{i}^{m}[j]${\tiny $_{L}$}$  = p_{i,i}^{\frac{m(m-1)}{2}}\sum \limits _{k=0}^{m} ( - 1)^{k} {\scriptsize\left(       \begin{array}{c} m\\
k\end{array}\right)}p_{i,j}^{k}p_{j,i}^{m-k} x_{i} ^{m - k}x_{j} x_{i} ^{k}$.

{\rm (iii)} $l_{i}^{m}[j]${\tiny $_{R}$}$  = p_{i,i}^{\frac{m(m-1)}{2}}\sum \limits _{k=0}^{m} ( - 1)^{k} {\scriptsize\left(       \begin{array}{c} m\\
k\end{array}\right)}p_{i,j}^{m-k}p_{j,i}^{k } x_{i} ^{m - k}x_{j} x_{i} ^{k}$.

{\rm (iv)}
$l_{i}^{m}[j]=\sum \limits _{k=0}^{m}(-1)^{k}p_{i,i}^{\frac{k(k-1)}{2}} p_{j,i}^{k}{\scriptsize\left(\begin{array}{cc} m\\
k \end{array}\right)_{p_{i,i}}}x_{i} ^{k }x_{j} x_{i} ^{m-k}$.

\end {Lemma}
\noindent {\it Proof.} It can be proved by induction. \qed

It is clear $l_{i}^{m}[j]${\tiny $_{L}$}$=l_{i}^{m}[j]${\tiny $_{R}$}$=p_{i,i}^{\frac{m(m-1)}{2}}p_{i,j}^{m}\overline{l}_{i}^{m}[j]^{-}$ if $p_{\cdot,\cdot}$ is a symmetrical bicharacter.

\begin {Lemma} \label {6.6} Assume that $\mathfrak B(V) $ is a  Nichols algebra of diagonal type and $m\in \mathbb N,i\neq j$.

Then {\rm (i)} If $p_{i,j}=p_{j,i}^{2}$ and $p_{j,i}=p_{i,j}^{2}$, i.e. $p_{i,j}p_{j,i}=1,p_{i,j}\in R_{3}$, then $l_{i}^{m}[j]${\tiny $_{R}$} $=0$.

{\rm (ii)} (1) $<y_{j},  l_{i}^{m}[j]${\tiny $_{R}$}$>= p_{i,i}^{\frac{m(m-1)}{2}}(p_{i,j}p_{j,i}^{-1}-p_{j,i})^{m}x_{i}^{m}$.

\ \ \ \ \ (2) $<y_{i}^{m}y_{j}, l_{i}^{m}[j]${\tiny $_{R}$}$ >=(p_{i,j}p_{j,i}^{-1}-p_{j,i})^{m} (m)_{p_{i,i}}^{!}$.

{\rm (iii)} (1)$<y_{i}^{k}, l_{i}^{m}[j]${\tiny $_{R}$}$> = p_{i,i}^{m-1}\{(p_{i,j}-p_{j,i} p_{i,j}^{-1}p_{i,i}^{k-m})(k)_{p_{i,i}^{-1}}<y_{i}^{k-1}, l_{i}^{m-1}[j]${\tiny $_{R}$}$>$

\hfill $+p_{i,j}p_{i,i}^{-k}x_{i}<y_{i}^{k}, l_{i}^{m-1}[j]${\tiny $_{R}$}$>
-p_{j,i}<y_{i}^{k}, l_{i}^{m-1}[j]${\tiny $_{R}$}$>x_{i}\}$ for $1\leq k \le m$.

(2)$<y_{i}^{m}, l_{i}^{m}[j]${\tiny $_{R}$}$>=(p_{i,j}-p_{j,i} p_{i,j}^{-1})^{m}(m)_{p_{i,i}}^{!} x_{j}$.

(3)$<y_{j}y_{i}^{m}, l_{i}^{m}[j]${\tiny $_{R}$}$>=(p_{i,j}-p_{j,i} p_{i,j}^{-1})^{m}(m)_{p_{i,i}}^{!}$.

{\rm (iv)} Assume that $p_{i,i}=1$. Then $l_{i}^{m}[j]${\tiny $_{R}$}$\not=0$,  when $p_{i,j}\neq p_{j,i}^{2}$ or $p_{j,i}\neq p_{i,j}^{2}$.

{\rm (v)} Assume that $p_{i,i}\not=1$. Then $l_{i}^{m}[j]${\tiny $_{R}$}$\not=0$,  when $ {\rm ord }(p_{i,i}) > m$ with $p_{i,j}\neq p_{j,i}^{2}$ or $p_{j,i}\neq p_{i,j}^{2}$.

\end {Lemma}
\noindent {\it Proof.}
{\rm (ii)} (1)$<y_{j},  l_{i}^{m}[j]${\tiny $_{R}$}
$ >= <y_{j},  p_{i,i}^{\frac{m(m-1)}{2}}\sum \limits _{k=0}^{m} ( - 1)^{k} {\scriptsize\left(       \begin{array}{c} m\\
k\end{array}\right)}p_{i,j}^{m-k}p_{j,i}^{k } x_{i} ^{m - k}x_{j} x_{i} ^{k} >$

$=p_{i,i}^{\frac{m(m-1)}{2}}p_{i,j}^{m}p_{j,i}^{-m }\sum \limits _{k=0}^{m} ( - 1)^{k} {\scriptsize\left(       \begin{array}{c} m\\
k\end{array}\right)}p_{i,j}^{-k}p_{j,i}^{2k }x_{i} ^{m}$

$=p_{i,i}^{\frac{m(m-1)}{2}}p_{i,j}^{m}p_{ji}^{-m
}(1-p_{i,j}^{-1}p_{j,i}^{2})^{m}x_{i} ^{m}$

$=p_{i,i}^{\frac{m(m-1)}{2}}(p_{i,j}p_{j,i}^{-1}-p_{j,i})^{m}x_{i}^{m}$.

(2) It can be proved by induction.

{\rm (iii)} We  show (1) by induction on $k$ for $1\leq k \le m$. See

$ <y_{i}, l_{i}^{m}[j]${\tiny $_{R}$} $> = <y_{i}, p_{i,i}^{m-1}p_{i,j} x_{i}l_{i}^{m-1}[j]${\tiny $_{R}$}$ - p_{i,i}^{m-1}p_{j,i} l_{i}^{m-1}[j]${\tiny $_{R}$}$x_i > $

$= p_{i,i}^{m-1}\{(p_{i,j}-p_{j,i} p_{i,j}^{-1}p_{i,i}^{1-m})l_{i}^{m-1}[j]${\tiny $_{R}$}$ +p_{i,j} p_{i,i}^{-1} x_i  <y_{i},  l_{i}^{m-1}[j]${\tiny $_{R}$}$> - p_{j,i} < y_i,  l_{i}^{m-1}[j]${\tiny $_{R}$}$>x_i  \}$.

Thus equation (1) holds when $k=1$. Assume $k>1$. See

$ <y_{i} ^{k}, l_{i}^{m}[j]${\tiny $_{R}$}$>$

$=<y_{i},  p_{i,i}^{m-1}\{(p_{i,j}-p_{j,i} p_{i,j}^{-1}p_{i,i}^{k-1-m})(k-1)_{p_{i,i}^{-1}}<y_{i}^{k-2}, l_{i}^{m-1}[j]${\tiny $_{R}$}$>+p_{i,j}p_{i,i}^{-k+1}x_{i}$

\hfill $<y_{i}^{k-1}, l_{i}^{m-1}[j]${\tiny $_{R}$}$>-p_{j,i}<y_{i}^{k-1}, l_{i}^{m-1}[j]${\tiny $_{R}$}$>x_{i}\} > $
(by inductive hypothesis)

$=p_{i,i}^{m-1}\{(p_{i,j}-p_{j,i} p_{i,j}^{-1}p_{i,i}^{k-1-m})(k-1)_{p_{i,i}^{-1}}<y_{i}^{k-1}, l_{i}^{m-1}[j]${\tiny $_{R}$}$>+p_{i,j}p_{i,i}^{-k+1}<y_{i}^{k-1}, l_{i}^{m-1}[j]${\tiny $_{R}$}$>+p_{i,j}p_{i,i}^{-k}x_{i}<y_{i}^{k}, l_{i}^{m-1}[j]${\tiny $_{R}$}$>-p_{j,i}<y_{i}^{k}, l_{i}^{m-1}[j]${\tiny $_{R}$}$>x_{i} -p_{j,i}p_{i,j}^{-1}p_{i,i}^{k-m}<y_{i}^{k-1}, l_{i}^{m-1}[j]${\tiny $_{R}$}$>\}$

$=$ the right hand side of (1). Consequently,  (1) holds.

Now we show (2) by induction on $m$.  (2) is clear when $m=1$.
By (1),  one obtains

$<y_{i}^{m}, l_{i}^{m}[j]${\tiny $_{R}$}$> =p_{i,i}^{m-1}(p_{i,j}-p_{j,i} p_{i,j}^{-1})(m)_{p_{i,i}^{-1}}<y_{i}^{m-1}, l_{i}^{m-1}[j]${\tiny $_{R}$}$>$

$=(p_{i,j}-p_{j,i} p_{i,j}^{-1})(m)_{p_{i,i}}<y_{i}^{m-1}, l_{i}^{m-1}[j]${\tiny $_{R}$}$>$

$=(p_{i,j}-p_{j,i} p_{i,j}^{-1})(m)_{p_{i,i}}(p_{i,j}-p_{j,i} p_{i,j}^{-1})^{m-1}(m-1)_{p_{i,i}}^{!}x_{j}$(by inductive hypothesis)

$=$the right hand side of  (2). Therefore,  (2) and (3) hold.

{\rm (iv)} If $p_{i,i}=1$,  then  $<y_{i}^{m}y_{j}, l_{i}^{m}[j]${\tiny $_{L}$}$>=(p_{i,j}p_{j,i}^{-1}-p_{j,i})^{m} (m)!$ and $<y_{j}y_{i}^{m}, l_{i}^{m}[j]${\tiny $_{L}$}$>=(p_{i,j}-p_{j,i} p_{i,j}^{-1})^{m}(m)!$ by {\rm (ii)} and {\rm (iii)}.

{\rm (v)}  It follows from  {\rm (ii)} and {\rm (iii)}. \qed

\begin {Lemma} \label {6.7} Assume that $\mathfrak B(V) $ is a  Nichols algebra of diagonal type and $m\in \mathbb N,i\neq j$.

Then {\rm (i)} If $p_{i,j}=p_{j,i}=1$, then $l_{i}^{m}[j]${\tiny $_{L}$} $=0$.

{\rm (ii)} (1) $<y_{j},  l_{i}^{m}[j]${\tiny $_{L}$}$>= p_{i,i}^{\frac{m(m-1)}{2}}  (1 - p_{i,j})^{m} x_{i} ^{m}$.
(2) $<y_{i}^{m}y_{j}, l_{i}^{m}[j]${\tiny $_{L}$}$ >=(1 - p_{i,j})^{m} (m)_{p_{i,i}}^{!}$.

{\rm (iii)} (1)$<y_{i}^{k}, l_{i}^{m}[j]${\tiny $_{L}$}$> = p_{i,i}^{m-1}\{(p_{j,i}-p_{i,i}^{k-m})(k)_{p_{i,i}^{-1}}<y_{i}^{k-1}, l_{i}^{m-1}[j]${\tiny $_{L}$}$>$

\hfill $+p_{j,i}p_{i,i}^{-k}x_{i}<y_{i}^{k}, l_{i}^{m-1}[j]${\tiny $_{L}$}$>
-p_{i,j}<y_{i}^{k}, l_{i}^{m-1}[j]${\tiny $_{L}$}$>x_{i}\}$ for $1\leq k \le m$.

(2)$<y_{i}^{m}, l_{i}^{m}[j]${\tiny $_{L}$}$>=(p_{j,i}-1)^{m}(m)_{p_{i,i}}^{!} x_{j}$.
(3)$<y_{j}y_{i}^{m}, l_{i}^{m}[j]${\tiny $_{L}$}$>=(p_{j,i}-1)^{m}(m)_{p_{i,i}}^{!}$.

{\rm (iv)} Assume that $p_{i,i}=1$. Then $l_{i}^{m}[j]${\tiny $_{L}$}$\not=0$,  when $p_{i,j}\not=1$ or $p_{j,i}\not=1$.

{\rm (v)} Assume that $p_{i,i}\not=1$. Then $l_{i}^{m}[j]${\tiny $_{L}$}$\not=0$,  when $ {\rm ord }(p_{i,i}) > m$ with $p_{i,j}\not=1$ or $ p_{j,i}\not=1$.
\end {Lemma}
\noindent {\it Proof.}
{\rm (ii)} (1)$<y_{j},  l_{i}^{m}[j]${\tiny $_{L}$}
$ >= <y_{j},  p_{i,i}^{\frac{m(m-1)}{2}}\sum \limits _{k=0}^{m} ( - 1)^{k} {\scriptsize\left(       \begin{array}{c} m\\
k\end{array}\right)}p_{i,j}^{k}p_{j,i}^{m-k} x_{i} ^{m - k}x_{j} x_{i} ^{k} >$

$= p_{i,i}^{\frac{m(m-1)}{2}}\sum \limits _{k=0}^{m} ( - 1)^{k} {\scriptsize\left(       \begin{array}{c} m\\
k\end{array}\right)}p_{i,j}^{k}p_{j,i}^{m-k}p_{ji}^{k-m} x_{i} ^{m } $

$=p_{i,i}^{\frac{m(m-1)}{2}}  (1 - p_{i,j})^{m} x_{i} ^{m}$.

(2) It can be proved by induction.

{\rm (iii)} We  show (1) by induction on $k$ for $1\leq k \le m$. See

$ <y_{i}, l_{i}^{m}[j]${\tiny $_{L}$} $> = <y_{i}, p_{i,i}^{m-1}p_{j,i} x_{i}l_{i}^{m-1}[j]${\tiny $_{L}$}$ - p_{i,i}^{m-1}p_{i,j} l_{i}^{m-1}[j]${\tiny $_{L}$}$x_i > $

$= p_{i,i}^{m-1}\{(p_{j,i}-p_{i,i}^{1-m})l_{i}^{m-1}[j]${\tiny $_{L}$}$ +p_{j,i} p_{i,i}^{-1} x_i  <y_{i},  l_{i}^{m-1}[j]${\tiny $_{L}$}$> - p_{i,j} < y_i,  l_{i}^{m-1}[j]${\tiny $_{L}$}$>x_i  \}$.

Thus equation (1) holds when $k=1$. Assume $k>1$. See

$ <y_{i} ^{k}, l_{i}^{m}[j]${\tiny $_{L}$}$>$

$=<y_{i},  p_{i,i}^{m-1}\{(p_{j,i}-p_{i,i}^{k-1-m})(k-1)_{p_{i,i}^{-1}}<y_{i}^{k-2}, l_{i}^{m-1}[j]${\tiny $_{L}$}$>+p_{j,i}p_{i,i}^{-k+1}x_{i}$

\hfill $<y_{i}^{k-1}, l_{i}^{m-1}[j]${\tiny $_{L}$}$>-p_{i,j}<y_{i}^{k-1}, l_{i}^{m-1}[j]${\tiny $_{L}$}$>x_{i}\} > $
(by inductive hypothesis)

$=p_{i,i}^{m-1}\{(p_{ji}-p_{i,i}^{k-1-m})(k-1)_{p_{i,i}^{-1}}<y_{i}^{k-1}, l_{i}^{m-1}[j]${\tiny $_{L}$}$>+p_{j,i}p_{i,i}^{-k+1}<y_{i}^{k-1}, l_{i}^{m-1}[j]${\tiny $_{L}$}$>+p_{j,i}p_{i,i}^{-k}x_{i}<y_{i}^{k}, l_{i}^{m-1}[j]${\tiny $_{L}$}$>-p_{i,j}<y_{i}^{k}, l_{i}^{m-1}[j]${\tiny $_{L}$}$>x_{i} -p_{i,j}p_{i,j}^{-1}p_{i,i}^{k-m}<y_{i}^{k-1}, l_{i}^{m-1}[j]${\tiny $_{L}$}$>\}$

$=$ the right hand side of (1). Consequently,  (1) holds.

Now we show (2) by induction on $m$.  (2) is clear when $m=1.$
By (1),  one obtains

$<y_{i}^{m}, l_{i}^{m}[j]${\tiny $_{L}$}$> = p_{i,i}^{m-1}(p_{j,i}-1)(m)_{p_{i,i}^{-1}}<y_{i}^{m-1}, l_{i}^{m-1}[j]${\tiny $_{L}$}$>$

$=(p_{j,i}-1)(m)_{p_{i,i}}<y_{i}^{m-1}, l_{i}^{m-1}[j]${\tiny $_{L}$}$>$

$=(p_{j,i}-1)(m)_{p_{i,i}}(p_{j,i}-1)^{m-1}(m-1)_{p_{i,i}}^{!} x_{j}$(by inductive hypothesis)

$=$the right hand side of  (2). Therefore,  (2) and (3) hold.

{\rm (iv)} If $p_{i,i}=1$,  then  $<y_{i}^{m}y_{j}, l_{i}^{m}[j]${\tiny $_{L}$}$>=(1 - p_{i,j})^{m} (m)!$ and $<y_{j}y_{i}^{m}, l_{i}^{m}[j]${\tiny $_{L}$}$>=(p_{j,i}-1)^{m}(m)!$ by {\rm (ii)} and {\rm (iii)}.

{\rm (v)}  It follows from  {\rm (ii)} and {\rm (iii)}. \qed

\begin {Proposition} \label {6.8} Assume that $\mathfrak B(V) $ is a  Nichols algebra of diagonal type. Then

{\rm (i)} $\dim (\mathfrak L(V)${\tiny $_{R}$}$)=\infty$ when there exist $i$ and $j$ with $i\not=j$, $p_{i,j}\neq p_{j,i}^{2}$ or $p_{j,i}\neq p_{i,j}^{2}$ and ${\rm ord}(p_{i,i})=1$ or $\infty$.

{\rm (ii)} $\dim (\mathfrak L(V)${\tiny $_{L}$}$)=\infty$ when there exist $i$ and $j$ with $i\not=j$, $p_{i,j}\not=1$ or $p_{j,i}\not=1$ and ${\rm ord}(p_{i,i})=1$ or $\infty$.
\end {Proposition}

\begin {Lemma} \label {6.9}
{\rm (See \cite [Lemma 5.1]{WZZ})}  The set of monomials $x_{i}^{s}x_{j}x_{i}^{t}$ with $ 0\leq s\leq m_{ij},0\leq t\leq{\rm ord }(p_{i,i})-1$ for $i<j$ and the set of monomials $x_{i}^{s}x_{j}x_{i}^{t}$ with $ 0\leq t\leq m_{ij},0\leq s\leq{\rm ord }(p_{i,i})-1$ for $i>j$ are linearly independent.
\end {Lemma}

We know $m_{ij}\leq {\rm ord }(p_{ii})-1$.

\begin {Proposition} \label {6.10} For $m\in \mathbb N,i\neq j$.

{\rm (i)} Assume that $p_{i,j}=p_{j,i}^{-1}\not=1$, then $l_{i}^{m}[j]${\tiny $_{L}$} $\neq0$
if and only if $ {\rm ord }(p_{i,i}) > m\geq0$, $l_{i}^{m}[j]${\tiny $_{R}$} $\neq0$
if and only if $p_{j,i}\notin R_{3}$ and $ {\rm ord }(p_{i,i}) > m\geq0$.

{\rm (ii)} Assume that $m_{ij}={\rm ord }(p_{i,i})-1,p_{i,j}p_{j,i}\neq1$, then
$l_{i}^{m}[j]${\tiny $_{L}$}$\not=0$ if and only if $0\leq m\leq m_{ij}-1+{\rm ord }(p_{i,i})$ if and only if $l_{i}^{m}[j]${\tiny $_{R}$}$\not=0$.
\end {Proposition}
\noindent {\it Proof.}
{\rm (i)}  We can show $l_{i}^{m}[j]${\tiny $_{L}$}$=p_{i,i}^{\frac{m(m-1)}{2}}(1-p_{i,j})^{m}ji^{m}$ and $l_{i}^{m}[j]${\tiny $_{R}$}$=p_{i,i}^{\frac{m(m-1)}{2}}p_{i,j}^{2m}(1-p_{j,i}^{3})^{m}ji^{m}$
by induction on $m$ since $ij=p_{i,j}ji$. the other is clear by \cite [Lemma 1.3.3]{He05}.

{\rm (ii)} We obtain $l_{i}^{m}[j]${\tiny $_{L}$}$  = p_{i,i}^{\frac{m(m-1)}{2}}\sum \limits _{k=0}^{m} ( - 1)^{k} {\scriptsize\left(       \begin{array}{c} m\\
k\end{array}\right)}p_{i,j}^{k}p_{j,i}^{m-k} x_{i} ^{m - k}x_{j} x_{i} ^{k}$ by Lemma \ref {6.5}{\rm (ii)}.
Then $l_{i}^{m_{ij}-1+{\rm ord }(p_{i,i})}[j]${\tiny $_{L}$}
$=l_{i}^{2{\rm ord }(p_{i,i})-2}[j]${\tiny $_{L}$}

$= p_{i,i}^{\frac{(2{\rm ord }(p_{i,i})-2)(2{\rm ord }(p_{i,i})-3)}{2}}\sum \limits _{k=0}^{2{\rm ord }(p_{i,i})-2} ( - 1)^{k} {\scriptsize\left(       \begin{array}{c} 2{\rm ord }(p_{i,i})-2\\
k\end{array}\right)}p_{i,j}^{k}p_{j,i}^{2{\rm ord }(p_{i,i})-2-k} x_{i} ^{2{\rm ord }(p_{i,i})-2 - k}x_{j} x_{i} ^{k}$

$= p_{i,i}^{\frac{(2{\rm ord }(p_{i,i})-2)(2{\rm ord }(p_{i,i})-3)}{2}} ( - 1)^{{\rm ord }(p_{i,i})-1} {\scriptsize\left(       \begin{array}{c} 2{\rm ord }(p_{i,i})-2\\
{\rm ord }(p_{i,i})-1\end{array}\right)}(p_{i,j}p_{j,i})^{{\rm ord }(p_{i,i})-1} x_{i} ^{{\rm ord }(p_{i,i})-1}x_{j} x_{i} ^{{\rm ord }(p_{i,i})-1}$

\noindent by \cite [Lemm 1.3.3(i)]{He05}. It is clear $l_{i}^{m_{ij}-1+{\rm ord }(p_{i,i})}[j]${\tiny $_{L}$}
$\neq 0$ since $x_{i} ^{{\rm ord }(p_{i,i})-1}x_{j} x_{i} ^{{\rm ord }(p_{i,i})-1}=x_{i} ^{m_{ij}}x_{j} x_{i} ^{{\rm ord }(p_{i,i})-1}$ or $x_{i} ^{{\rm ord }(p_{i,i})-1}x_{j}x_{i} ^{m_{ij}}$ is a basic element by  Lemma \ref {6.5}. Then $l_{i}^{m}[j]${\tiny $_{L}$}$\neq 0$ if $0\leq m\leq m_{ij}-1+{\rm ord }(p_{i,i})$. On
the other hand,

\noindent $l_{i}^{m_{ij}+{\rm ord }(p_{i,i})}[j]${\tiny $_{L}$}

\noindent $=
p_{i,i}^{\frac{(2{\rm ord }(p_{i,i})-2)(2{\rm ord }(p_{i,i})-3)}{2}} ( - 1)^{{\rm ord }(p_{i,i})-1} {\scriptsize\left(       \begin{array}{c} 2{\rm ord }(p_{i,i})-2\\
{\rm ord }(p_{i,i})-1\end{array}\right)}(p_{i,j}p_{j,i})^{{\rm ord }(p_{i,i})-1} [x_{i},x_{i} ^{{\rm ord }(p_{i,i})-1}x_{j} x_{i} ^{{\rm ord }(p_{i,i})-1}]${\tiny $_{L}$}

\noindent $=0$. The rest of the proof is similar. \hfill $\Box$

\begin {Lemma} \label {6.11} Assume that $\mathfrak B(V)$ is a connected Nichols algebra of diagonal type with $\dim V>1$. Then

{\rm (i)} If $p _{i,i}  = p_{j,j} =-1$ and ${\rm ord } (p_{i,j} p _{j,i}) = \infty$ for $i\neq j$.  then $\dim (\mathfrak L(V)${\tiny $_{R}$}$) = \infty$,
$\dim (\mathfrak L(V)${\tiny $_{L}$}$) = \infty$.

{\rm (ii)} If $\Delta (\mathfrak B(V))$ is an arithmetic root system and  there exists $u\in D$ such that ${\rm ord } (p_{u,u}) = \infty$, then there exists $1\le i \le n$ such that ${\rm ord } (p_{i,i}) = \infty$ or $\dim (\mathfrak L(V)${\tiny $_{R}$}$) = \infty$, $\dim (\mathfrak L(V)${\tiny $_{L}$}$) = \infty$.

\end {Lemma}
\noindent{\it Proof.}  { \rm (i)} It follows from Theorem \ref {4.3} and Theorem \ref {4.5}.

{\rm (ii)} If there exists $u \in D$  such that ${\rm ord }  (p_{u,u}) =\infty$, then there exists a $1\le i \le n$  such that ${\rm ord }  (p_{i,i}) = \infty $  but Case of  \cite  [Row 3 Diagram 2, Table 1] {He05},  \cite  [ Row 8, Diagram 2,  Table 2; Row 9, Diagram 4,  Table 2;  Row 10, Diagram 3,  Table 2] {He05} and   \cite  [ Row 10, Diagram 6,  Appendix B;  Row 12, Diagram 5,  Appendix B. Row 2,   Appendix C; Row 10, Diagram 2,   Appendix C  ] {He06a}.  By \rm {(i)}, $\dim (\mathfrak L(V)${\tiny $_{R}$}$) = \infty$, $\dim (\mathfrak L(V)${\tiny $_{L}$}$) = \infty$.  \qed

\subsection {Finiteness of $\mathfrak L(V)${\tiny $_{L}$} and $\mathfrak L(V)${\tiny $_{R}$}}

\begin {Proposition} \label {6.13} Assume that $\mathfrak B(V) $ is a connected  Nichols algebra of diagonal type. Then $\dim(\mathfrak B(V))=\infty$ and  $ \dim(\mathfrak L(V)${\tiny $_{L}$}$)<\infty$  if and only if   $\dim V=1$ and ${\rm ord}(p_{1,1})=1$ or $\infty$.
\end {Proposition}

\noindent{\it Proof.}
Sufficiency. It is clear.

Necessity. Assume $\dim V>1.$ By Theorem \ref {6.4},  $ \Delta (\mathfrak B(V))$ is an arithmetic root system.  Consequently,  there exists $u \in D$ such that ${\rm ord } (p_{u,u}) = \infty$. Considering Lemma \ref {6.11} (ii) and Proposition \ref {6.8} (ii) we get a contradiction. \qed

\begin {Proposition} \label {6.14} Assume that $\mathfrak B(V) $ is a connected  Nichols algebra of diagonal type. Then $\dim(\mathfrak B(V))=\infty$ and $\dim(\mathfrak L(V)${\tiny $_{R}$}$)<\infty$  if and only if $\dim V=1,{\rm ord}(p_{1,1})=1$ or $\infty$.
\end {Proposition}
\noindent {\it Proof.}   Sufficiency. It is clear. Assume $\dim V>1.$ By Theorem \ref {6.4},  $ \Delta (\mathfrak B(V))$ is an arithmetic root system.  Consequently,  there exists $u \in D$ such that ${\rm ord } (p_{u,u}) = \infty$. Considering Lemma \ref {6.11} (ii) and Proposition \ref {6.8} (i) we have $p_{i,j}= p_{j,i}^{2},p_{j,i}= p_{i,j}^{2}$, then $p_{i,j}p_{j,i}=1,p_{i,j}\in R_{3}$, which is a contradiction. \qed

Using Proposition \ref {6.13} and  Proposition \ref {6.14}, we obtain the following result.

\begin {Theorem} \label {6.15} Assume that $\mathfrak B(V)$ is a connected Nichols algebra of diagonal type with $\dim V>1$, then the following conditions are equivalent:
{\rm (i)} $\mathfrak B(V)$ is finite-dimensional;
{\rm (ii)} $\mathfrak L(V)${\tiny $_{L}$} is finite-dimensional;
{\rm (iii)} $\mathfrak L(V)${\tiny $_{R}$} is finite-dimensional.
\end {Theorem}

\textbf{Weicai Wu}\\
School of Mathematics and Statistics, Guangxi Normal University,
541004 Guilin, Guangxi, People's Republic of China.\\
E-mail: weicaiwu@hnu.edu.cn

\end{document}